\documentclass[12pt,reqno]{amsart}
\usepackage{amsmath,amssymb}
\begin{document}
 \baselineskip=17pt  

\title{Random iteration with place dependent probabilities}
\author{Rafa\l\ Kapica }
\address{AGH University of Science and Technology, Faculty of Applied Mathematics,
al. A. Mickiewicza 30, 30-059 Krakow, Poland}
\email{rafal.kapica@agh.edu.pl}

\author{Maciej \' Sl\c eczka }
\address{Institute of Mathematics, University of Silesia in Katowice, Bankowa 14,
40-007 Katowice, Poland}
\email{sleczka@math.us.edu.pl}


\renewcommand{\thefootnote}{}



\renewcommand{\thefootnote}{\arabic{footnote}}
\setcounter{footnote}{0}



\renewcommand*{\abstractname}{Abstract}
\begin{abstract}
Markov chains arising from random iteration of functions $S_{\theta}:X\to X$, $\theta \in \Theta$, where $X$ is a Polish space and $\Theta$ is an arbitrary set of indices are considerd. At $x\in X$, $\theta$ is sampled from a distribution $\vartheta_x$ on $\Theta$ and $\vartheta_x$ are different for different $x$. Exponential convergence to a unique invariant measure is proved. This result is applied to the case of random affine transformations on ${\mathbb R}^d$ giving the existence of exponentially attractive perpetuities with place dependent probabilities.
\end{abstract}

\maketitle

\newtheorem{theorem}{Theorem}[section]
\newtheorem{lemma}{Lemma}[section]
\newtheorem{proposition}{Proposition}[section]
\newtheorem{remark}{Remark}[section]
\newtheorem{corollary}{Corollary}[section]
\newtheorem{definition}{Definition}[section]

\section{Introduction}

We consider the Markov chain of the form $ X_0=x_0$, $X_1=S_{\theta_0}(x_0)$, $X_2=S_{\theta_1}\circ S_{\theta_0}(x_0)$ and inductively
\begin{equation}\label{ex}
X_{n+1}=S_{\theta_n}(X_n),
\end{equation}
where $S_{\theta_0}$, $S_{\theta_1}$,...,$S_{\theta_n}$ are randomly chosen from a family $\{S_{\theta}:\, \theta\in\Theta\}$ of functions that map a state space $X$ into itself. If the chain is at $x\in X$ then $\theta\in\Theta$ is sampled from a distribution $\vartheta_x$ on $\Theta$, where $\vartheta_x$ are, in general, different for different $x$. We are interested in the rate of convergence to a stationary distribution $\mu_*$ on $X$, i.e.
\begin{equation}\label{convergence}
P\{X_n\in A\}\to \mu_*(A)\qquad\text{as}\qquad n\to\infty.
\end{equation}
In the case of constant probabilities, i.e. $\vartheta_x=\vartheta_y$ for $x,y\in X$, the basic tool when studying asymptotics of (\ref{ex}) are backward iterations 
$$
Y_{n+1}=S_{\theta_0}\circ S_{\theta_1}\circ ...\circ S_{\theta_n}(x).
$$
Since $X_n$ and $Y_n$ are identically distributed and, under suitable conditions, $Y_n$ converge almost surely at exponential rate to some random element $Y$, one obtains exponential convergence in (\ref{convergence}) (see \cite{DF} for bibliography and excellent survey of the field). For place dependent $\vartheta_x$ we need a different approach because distributions of $X_n$ and $Y_n$ are not equal.
\newline
The simplest case when $\Theta=\{1,...,n\}$ is treated in \cite{BDEG} and \cite{Sz}, where the existence of a unique attractive invariant measure is established. Similar result holds true when $\Theta=[0,T]$ and $\vartheta_x$ are absolutely continuous (see \cite{HSz}). Recently it was shown that the rate of convergence in the case of $\Theta=\{1,...,n\}$ is exponential (see \cite{S}).
\newline
In this paper we treat the general case of place dependent $\vartheta_x$ for arbitrary $\Theta$ and prove the existence of a unique exponentially attractive invariant measure for (\ref{ex}). Our approach is based on the coupling method which can be briefly described as follows. For arbitrary starting points $x,{\bar x}\in X$ we consider chains $(X_n)_{n\in\mathbb{N}_0}$, $({\bar X}_n)_{n\in\mathbb{N}_0}$ with $X_0=x_0$, ${\bar X}_0={\bar x}_0$ and try to build correlations between $(X_n)_{n\in\mathbb{N}_0}$ and $({\bar X}_n)_{n\in\mathbb{N}_0}$ in order to make their trajectories as close as possible. This can be done because the transition probability function ${\mathbf B}_{x,y}(A)=P\{ (X_{n+1},{\bar X}_{n+1})\in A\, |\, (X_{n},{\bar X}_{n})=(x,y)\}$ of the coupled chain $(X_n,{\bar X}_n)_{n\in{\mathbb N}_0}$ taking values in $X^2$ can be decomposed (see \cite{H}) in the following way
$$
{\mathbf B}_{x,y}={\mathbf Q}_{x,y}+{\mathbf R}_{x,y},
$$ 
where subprobability measures ${\mathbf Q}_{x,y}$ are contractive in metric $d$ on $X$:
$$
\int_{X^2} d(u,v)\, {\mathbf Q}_{x,y}(du,dv)\le\alpha d(x,y)
$$
for some constant $\alpha\in (0,1)$. 
\newline
Since transition probabilities for (\ref{ex}) can be mutually singular for even very close points, one cannot expect that chains $(X_n)_{n\in\mathbb{N}_0}$ and $({\bar X}_n)_{n\in\mathbb{N}_0}$ couple in finite time ($X_n={\bar X}_n$ for some $n\in {\mathbb N}_0$) as in classical coupling constructions (\cite{Lindvall}) leading to the convergence in the total variation norm. On the contrary, they only couple at infinity ($d(X_n,{\bar X}_n)\to 0$ as $ n\to \infty$) so this method is sometimes called asymptotic coupling (\cite{HMSch}) and gives convergence in *-weak topology.
\newline  
The paper is organized as follows. In Section 2 we formulate and prove theorem which assures exponential convergence to an invariant measure for a class of Markov chains. This theorem is applied in Section 3 to chains generated by random iteration of functions. In Section 4 we discuss special class of such functions, random affine transformations on ${\mathbb R}^d$, thus generalizing the notion of {\it perpetuity} to the place dependent case.

\section{An exponential convergence result}

\subsection{Notation and basic definitions.}

Let $(X,d)$ be a \emph {Polish} space, i.e. a complete and separable metric space and denote by $\mathcal{B}_X$ the 
$\sigma$-algebra of Borel subsets of $X$. By $B_b(X)$ we denote the space of bounded Borel-measurable functions equipped with the supremum norm, $C_b(X)$ stands for the subspace of bounded continuous functions. Let 
$\mathcal{M}_{fin}(X)$ and $\mathcal{M}_1(X)$ be the sets of Borel measures on $X$ such that $\mu(X)<\infty$ for $\mu\in\mathcal{M}_{fin}(X)$ and $\mu(X)=1$ for $\mu\in\mathcal{M}_1(X)$. The elements of $\mathcal{M}_1(X)$ are called \emph {probability measures}. The elements of $\mathcal{M}_{fin}(X)$ for which  $\mu(X)\le 1$ are called \emph {subprobability measures}. By $supp\, \mu$ we denote the support of the measure $\mu$. We also define
$$
{\mathcal M}_1^L (X)=\{\mu\in{\mathcal M}_1(X):\,\int_X L(x)\mu (d x)<\infty\}
$$
where $L:X\to [0,\infty)$ is an arbitrary Borel measurable function and
$$
{\mathcal M}_1^1 (X)=\{\mu\in{\mathcal M}_1(X):\,\int_X d(\bar x ,x)\mu (d x)<\infty\},
$$
where $\bar x\in X$ is fixed. By the triangle inequality the definition of ${\mathcal M}_1^1 (X)$ is independent of the choice of $\bar x$.
\newline
The space $\mathcal{M}_1(X)$ is equipped with the \emph{Fortet-Mourier metric}:
$$
\Vert\mu_1-\mu_2\Vert_{FM}=\sup\{|\int_X f(x)(\mu_1-\mu_2)(dx)|:\, f\in\mathcal{F}\},
$$
where
$$
\mathcal{F}=\{f\in C_b(X):\, |f(x)-f(y)|\le 1 \quad\text{and}\quad |f(x)|\le 1\quad\text{for}\quad x,y\in X\}.
$$
The space $\mathcal{M}_1^1(X)$ is equipped with the \emph{Wasserstein metric}:
$$
\Vert\mu_1-\mu_2\Vert_{W}=\sup\{|\int_X f(x)(\mu_1-\mu_2)(dx)|:\, f\in\mathcal{W}\},
$$
where
$$
\mathcal{W}=\{f\in C_b(X):\, |f(x)-f(y)|\le 1 \quad\text{for}\quad x,y\in X\}.
$$
By $\Vert\cdot\Vert$ we denote the total variation norm. If a measure $\mu$ is nonnegative then $\Vert \mu\Vert$ is simply the total mass of $\mu$.
\newline
Let $P:B_b(X)\to B_b(X)$ be a \emph{Markov operator}, i.e. a linear operator satisfying $P{\bf 1}_X={\bf 1}_X$ and $Pf(x)\ge 0$ if $f\ge 0$. Denote by $P^{*}$ the the dual operator, i.e operator $P^{*}:\mathcal{M}_{fin}(X)\to
\mathcal{M}_{fin}(X)$ defined as follows 
$$
P^{*}\mu(A):=\int_X P {\bf 1}_A(x)\mu(dx)\qquad\text{for}\qquad A\in\mathcal{B}_X.
$$
We say that a measure $\mu_*\in\mathcal{M}_1(X)$ is \emph{invariant} for $P$ if
$$
\int_X Pf(x)\mu_*(dx)=\int_X f(x)\mu_*(dx)\qquad\text{for every}\qquad f\in B_b(X)
$$
or, alternatively, we have $P^* \mu_*=\mu_*$. \newline
By $\{\mathbf{P}_x:\, x\in X\}$ we denote a \emph{transition probability function} for $P$, i.e. a family of measures $\mathbf{P}_x\in\mathcal{M}_1(X)$ for $x\in X$ such that the map $x\mapsto\mathbf{P}_x(A)$ is measurable for every $A\in\mathcal{B}_X$ and 
$$
Pf(x)=\int_X f(y) \mathbf{P}_x(dy)\qquad\text{for}\qquad x\in X\quad\text{and}\quad f\in B_b(X)
$$
or equivalently $P^*\mu(A)=\int_X \mathbf{P}_x(A)\mu(dx)$ for $A\in\mathcal{B}_X$ and $\mu\in\mathcal{M}_{fin}(X)$.

\subsection{Formulation of the theorem.}

\begin{definition}
A coupling for $\{\mathbf{P}_x: x\in X\}$ is a family $\{\mathbf{B}_{x,y}:\, x,y\in X\}$ of probability measures on $X\times X$ such that for every $B\in\mathcal{B}_{X^2}$ the map $X^2\ni (x,y)\mapsto \mathbf{B}_{x,y}(B)$ is measurable and
$$ 
\mathbf{B}_{x,y}(A\times X)=\mathbf{P}_x(A),\qquad \mathbf{B}_{x,y}(X\times A)=\mathbf{P}_y(A)
$$
for every $x,y\in X$ and $A\in\mathcal{B}_X$.
\end{definition}
In the following we assume (see \cite{H}) that there exists a family $\{\mathbf{Q}_{x,y}:\,x,y\in X\}$ of subprobability measures on $X^2$ such that the map $(x,y)\mapsto \mathbf{Q}_{x,y}(B)$ is measurable for every Borel $B\subset X^2$ and
$$
\mathbf{Q}_{x,y}(A\times X)\le \mathbf{P}_x(A)\qquad\text{and}\qquad \mathbf{Q}_{x,y}(X\times A)\le \mathbf{P}_y(A)
$$
for every $x,y\in X$ and Borel $A\subset X$.
\newline
Measures $\{\mathbf{Q}_{x,y}:x,y\in X\}$ allow us to construct a coupling for $\{\mathbf{P}_x:x\in X\}$. Define on $X^2$ the family of measures $\{\mathbf{R}_{x,y}:x,y\in X\}$ which on rectangles $A\times B$ are given by
$$
\mathbf{R}_{x,y}(A\times B)=\frac{1}{1-\mathbf{Q}_{x,y}(X^2)}(\mathbf{P}_x(A)-\mathbf{Q}_{x,y}(A\times X)) (\mathbf{P}_y(B)-\mathbf{Q}_{x,y}(X\times B)),
$$
when $\mathbf{Q}_{x,y}(X^2)<1$ and $\mathbf{R}_{x,y}(A\times B)=0$ otherwise. A simple computation shows that the family $\{\mathbf{B}_{x,y}:\, x,y\in X\}$ of measures on $X^2$ defined by
\begin{equation}\label{BQR}
\mathbf{B}_{x,y}=\mathbf{Q}_{x,y}+\mathbf{R}_{x,y}\quad\text{for}\quad x,y\in X 
\end{equation}
is a coupling for $\{\mathbf{P}_x:\, x\in X\}$.
\newline
For every $r>0$ define $D_r =\{(x,y)\in X^2:\,\,d(x,y)<r\,\,\}$.
\newline
Now we list our assumptions on Markov operator $P$ and transition probabilities $\{\mathbf{Q}_{x,y}:\,x,y\in X\}$.
\newline
{\bf A0} $P$ {\it is a} \emph{ Feller operator}, {\it i.e.} $P(C_b(X))\subset C_b (X)$. 
\newline
{\bf A1} {\it There exists a} \emph{ Lyapunov function}\, {\it for} $P$, {\it i.e. continuous function} $L: X\to [0, \infty )$     {\it such that} 
$L$ {\it is bounded on bounded sets, }
$\lim _{x\to \infty}L(x)=+\infty$ ({\it for bounded $X$ this condition is omitted}) {\it and for some} $\lambda \in (0,1), \, c>0$
$$
PL (x)\le \lambda L(x) +c \qquad for \qquad x\in X.
$$
{\bf A2} {\it There exist} $F\subset X^2$ {\it and} $\alpha \in  (0,1)$ {\it such that} $supp \,{\mathbf Q}_{x,y}\subset F$ {\it and}
\begin{equation}\label{A2}
\int_{ X^2} d (u,v){\mathbf Q}_{x,y}(d u,d v) \le \alpha d(x,y)\qquad for \qquad (x,y)\in F.
\end{equation}
{\bf A3} {\it There exist} $\delta>0,\, l>0$ {\it and} $\nu\in (0,1]$ {\it such that} 
\begin{equation}\label{A3}
1-\| {\mathbf Q}_{x,y} \| \le l d(x,y)^{\nu}\ \qquad\text{and}\qquad {\mathbf Q}_{x,y} (D_{\alpha d(x,y)} )\ge\delta
\end{equation}
{\it for} $(x,y)\in F$.
\newline
{\bf A4} {\it There exist} $\beta\in (0,1)$, ${\tilde C}>0$ {\it and} $R>0$ {\it such that for} 
$$
\kappa (\,(x_n,y_n)_{n\in{\mathbb N}_0}\,)=\inf \{n\in {\mathbb N}_0 :\, (x_n,y_n)\in F\quad\text{and}\quad L(x_n)+L(y_n)<R\}
$$
{\it we have}
$$
{\mathbb E}_{x,y} \beta^{-\kappa}\le {\tilde C}\qquad whenever \qquad L(x)+L(y)<\frac{4c}{1-\lambda},
$$
{\it where} $\mathbb{E}_{x,y}$ {\it denotes here the expectation with respect to the chain starting from} $(x,y)$ {\it and with trasition function } $\{\mathbf{B}_{x,y}:\, x,y\in X\}$.
\newline

{\it Remark.} Condition {\bf A4} means that the dynamics quickly enters the domain of contractivity $F$. In this paper we discuss Markov chains generated by random iteration of functions for which always $F=X^2$ and $L(x)=d(x,{\bar x})$ with some fixed ${\bar x}\in X$, so {\bf A4} is trivially fulfilled when $R=\frac{4c}{1-\lambda}$. There are, however, examples of random dynamical systems for which $F$ is a proper subset of $X^2$. Indeed, in {\it contractive Markov systems} introduced by I. Werner in \cite{W} we have $X=\sum_{i=1}^{n} X_i$ but $F=\sum_{i=1}^n X_i\times X_i$. They are studied in \cite{Exp}.\newline

Now we formulate the main result of this section. Its proof is given in Section 2.4.
\begin{theorem}\label{t1} Assume {\bf A0} -- {\bf A4}. Then operator $P$ possesses a unique invariant measure 
$\mu_{*}\in\mathcal{M}_1^L (X)$, which is attractive, i.e.
$$
\lim\limits_{n\to\infty}\int_X P^nf(x)\, \mu(dx)=\int_X f(x)\, \mu(dx)\quad\text{for}\quad f\in C_b(X),\, \mu\in\mathcal{M}_1(X).
$$
Moreover, there exist $q\in (0,1)$ and $C>0$ such that
\begin{equation}\label{mix}
\| P^{* n}\mu -\mu_{*}\|_{FM}\le q^n C(1+\int_X L(x)\mu (dx) )
\end{equation}
for $\mu\in\mathcal{M}_1^L (X)$ and $n\in\mathbb{N}$.
\end{theorem}

{\it Remark.} In \cite{HMSch}, Theorem 4.8, authors formulate sufficient conditions for the existence of a unique exponentially attractive invariant measure for continuous-time Markov semigroup $\{P(t)\}_{t\ge 0}$, that do not refer to coupling. One of assumptions is that there exists distance-like (i.e. symmetric, lower semi-continuous and vanishing only on the diagonal) function 
$d:X\times X\to [0,1]$ which is contractive for some $P(t_*)$, i.e. there exists $\alpha<1$ such that for every $x,y\in X$ with $d(x,y)<1$ we have
$$
d(\mathcal{P}(x,\cdot),\mathcal{P}(y,\cdot))\le\alpha d(x,y),
$$
where $\mathcal{P}(\cdot,\cdot):X\times\mathcal{B}_X \to [0,1]$ is transition kernel for $P(t_*)$. This assumption is stronger than 
{\bf A2}, since measures $\mathbf{R}_{x,y}$ in (\ref{BQR}) need not be contractive (i.e. 
$\int_{X^2} d(u,v)\, \mathbf{R}_{x,y} (du,dv)\le \alpha d(u,v)$) for any distance-like function $d$.\newline

\subsection{Measures on the pathspace.}

For fixed $(x_0,y_0)\in X^2$ the next step of a chain with transition probability function $\mathbf{B}_{x,y}=\mathbf{Q}_{x,y}+\mathbf{R}_{x,y}$ can be drawn according to $\mathbf{Q}_{x_0,y_0}$ or according to $\mathbf{R}_{x_0,y_0}$. To distinguish these two cases we introduce (see \cite{H}) the augmented space $\widehat{X}=X^2\times\{0,1\}$ and the transition function $\{\widehat{\mathbf{B}}_{x,y,\theta}:\, (x,y,\theta)\in \widehat{X}\}$ on $\widehat{X}$ given by
$$
\widehat{\mathbf{B}}_{x,y,\theta}= \widehat{\mathbf{Q}}_{x,y,\theta} + \widehat{\mathbf{R}}_{x,y,\theta},
$$
where $\widehat{\mathbf{Q}}_{x,y,\theta}=\mathbf{Q}_{x,y}\times\delta_1$ and $\widehat{\mathbf{R}}_{x,y,\theta}=\mathbf{R}_{x,y}\times\delta_0$.
The parameter $\theta\in\{0,1\}$ is responsible for choosing measures $\mathbf{Q}_{x,y}$ and $\mathbf{R}_{x,y}$. If a Markov chain with transition function $\{\widehat{\mathbf{B}}_{x,y,\theta}:\, (x,y,\theta)\in \widehat{X}\}$ stays in the set $X^2\times\{1\}$ at time $n$ it means that the last step was drawn according to $\mathbf{Q}_{u,v}$, for some $(u,v)\in X^2$.
\newline
For every $x\in X$ finite-dimensional distributions $\mathbf{P}_x^{0,...,n}\in\mathcal{M}_1(X^{n+1})$ are defined by
$$
\mathbf{P}_x^{0,...,n}(B)=\int_X\mu(d x_0)\int_X\mathbf{P}_{x_1}(d x_2)...\int_X\mathbf{P}_{x_{n-1}}(d x_n){\bf 1}_B(x_0,...,x_n)
$$
for $n\in\mathbb{N}_0$, $B\in\mathcal{B}_{X^{n+1}}$. By the Kolmogorov extension theorem we obtain the measure $\mathbf{P}_x^{\infty}$ on the pathspace $X^{\infty}$. Similarly we define measures $\mathbf{B}_{x,y}^{\infty}$, $\widehat{\mathbf{B}}_{x,y,\theta}^{\infty}$ on $(X\times X)^{\infty}$ and $\widehat{X}^{\infty}$. 
These measures have the following interpretation. Consider the Markov chain $(X_n,Y_n)_{n\in\mathbb{N}_0}$ on $X\times X$, starting from $(x_0,y_0)$, with the transition function $\{\mathbf{B}_{x,y}:\, x,y\in X\}$, obtained by canonical Kolmogorov construction, i.e. $\Omega=(X\times X)^{\infty}$ is the sample space equipped with the probability measure $\mathbb{P}=\mathbf{B}_{x_0,y_0}^{\infty}$, $X_n(\omega)=x_n$, $Y_n(\omega)=y_n$, where $\omega=(x_k,y_k)_{k\in\mathbb{N}_0}\in\Omega$, and $n\in\mathbb{N}_0$. 
Then $(X_n)_{n\in\mathbb{N}_0}$, $(Y_n)_{n\in\mathbb{N}_0}$ are Markov chains in $X$, starting from $x_0$ and $y_0$, with the transition function $\{\mathbf{P}_{x}:\, x\in X\}$, and $\mathbf{P}_x^{\infty}$, $\mathbf{P}_y^{\infty}$ are their measures on the pathspace $X^{\infty}$.
\newline 
In this paper we often consider marginals of measures on the pathspace. If $\mu$ is a measure on a measurable space $X$ and $f:X\to Y$ is a measurable map, then $f^{\#}\mu$ is the measure on $Y$ defined by $f^{\#}\mu (A)=\mu (f^{-1}(A))$. So, if we denote by $pr$ the projection map from  a product space to its component, then $pr^{\#}\mu$ is simply the marginal of $\mu$ on this component.
\newline
In the following we consider Markov chains 
on $\widehat{X}$ with the transition function $\{\widehat{\mathbf{B}}_{x,y,\theta}:\, x,y\in X, \theta\in\{0,1\}\}$. 
We adopt the convention that $\theta _0=1$, so $\Phi$ always starts from $X^2\times\{1\}$, and define 
$$
\widehat{\mathbf{B}}_{x,y}^{\infty}:=\widehat{\mathbf{B}}_{x,y,1}^{\infty}.
$$ 
\newline
For $b\in\mathcal{M}_{fin}(X^2)$ we write 
$$\widehat{\mathbf{B}}_b^{\infty}(B)=\int_X \widehat{\mathbf{B}}_{x,y}^{\infty}(B)\, b(dx,dy),\qquad B\in\mathcal{B}_{\widehat{X}^{\infty}},
$$ 
$$
\mathbf{Q}_b(A)=\int_{X^2}\mathbf{Q}_{x,y}(A)\, b(dx,dy),\qquad A\in\mathcal{B}_{X^2}
$$
and
$$
\mathbf{Q}^n_{x,y}(A)=\mathbf{Q}_{\mathbf{Q}^{n_0 - 1}_{x,y}} (A), \qquad A\in\mathcal{B}_{X^2}.
$$
When studying the asymptotics of the a chain $(X_n)_{n\in\mathbb{N}_0}$ with a transition function $\{\mathbf{P}_x:\, x\in X\}$ it is particularly interesting whether a coupled chain $(X_n,Y_n)_{n\in\mathbb{N}_0}$ is moving only according to the contractive part $\mathbf{Q}_{x,y}$ of the transition function $\mathbf{B}_{x,y}$. For every subprobability measure $b\in\mathcal{M}_{fin}(X^2)$ we define finite-dimensional subprobability distributions $\mathbf{Q}_b^{0,...,n}\in\mathcal{M}_{fin}((X\times X)^{n+1})$
$$
\begin{aligned}
\mathbf{Q}_b^{0,...,n}(B)=&\int_{X^2}b(dx_0,dy_0)\int_{X^2}\mathbf{Q}_{x_0,y_0}(dx_1,dy_1)...\\
&...\int_{X^2}\mathbf{Q}_{x_{n-1},y_{n-1}}{\bf1}_B((x_0,y_0),...,(x_n,y_n)),
\end{aligned}
$$
where $B\in\mathcal{B}_{(X\times X)^{n+1}}$, $n\in\mathbb{N}_0$. 
Since the family $\{\mathbf{Q}_b^{0,...,n}:\, n\in\mathbb{N}_0\}$ need not be consistent, we cannot use the Kolmogorov extension theorem to obtain a measure on the whole pathspace $\widehat{X}^{\infty}$. 
However, defining for every $b\in\mathcal{M}_{fin}(X^2)$ the measure $\mathbf{Q}_b^{\infty}\in\mathcal{M}_{fin}(\widehat{X}^{\infty})$ by
$$
\mathbf{Q}_b^{\infty}(B)=\widehat{\mathbf{B}}_b^{\infty}(B\cap(X^2\times \{1\})^{\infty}),
$$
where $B\in\mathcal{B}_{\widehat{X}^{\infty}}$, one can easily check that for every cylindrical set $B=A\times {\widehat X}^{\infty}$, $A\in\mathcal{B}_{\widehat{X}^n}$, we have
\begin{equation}\label{Qinfty}
\mathbf{Q}_b^{\infty}(B)=\lim_{n\to\infty }\mathbf{Q}_b^{0,...,n}(pr_{(X^2)^{n+1}}(A)).
\end{equation}

\subsection{Proof of Theorem \ref{t1}.}

Before proceeding to the proof of Theorem  \ref{t1} we formulate two lemmas. The first one is partially inspired by the reasoning which can be found in \cite{O}.

\begin{lemma} \label{ol} Let $Y$ be a metric space and let $(Y_n^y)_{n\in\mathbb{N}_0}$ be a family of Markov chains indexed by starting point $y\in Y$, with common transition function $\{\pi_y:\,y\in Y\}$. Let $V:Y\to [0,\infty)$ be a Lapunov function for $\{\pi_y:\,y\in Y\}$. 
Assume that for some bounded and measurable $A\subset Y$ there exist $\lambda\in (0,1)$ and $C_{\rho}>0$ such that for 
$$
\rho ((y_n)_{n\in{\mathbb N}_0})=\inf\{n\ge 1:\,y_n\in A\}
$$
we have
$$
\mathbb{E}_{y} \lambda^{-\rho}\le C_{\rho}(V(y_0)+1)\quad\text{for}\quad y\in Y,
$$
where $\mathbb{E}_y$ is the expectation with respect to the measure $\mathbb{P}_y$ on $Y^{\infty}$ induced by $(Y_n^y)_{n\in\mathbb{N}_0}$. Moreover, assume that for some measurable $B\subset Y$ and 
$$
\epsilon((y_n)_{n\in \mathbb{N}_0})=\inf\{n \ge 1 :\, y_n\notin B\}
$$
there exist constants $p>0$, $\beta\in (0,1)$ and $C_{\epsilon}>0$  such that 
$$
\mathbb{P}_y(\{(y_n)_{n\in\mathbb{N}_0}:\, \forall_{n\ge 1}\,y_n\in B\})>p
\qquad\text{and}\qquad
\mathbb{E}_y {\bf 1}_{\{\epsilon<\infty\}} \beta^{-\epsilon}\le C_{\epsilon},
$$
for every $y\in A$.\newline
Then there exist $\gamma \in (0,1)$ and $C>0$ such that for 
$$
\tau((y_n)_{n\in\mathbb{N}_0})=\inf\{n \ge 1 :\, \forall_{k\ge n} \,y_k\in B\}
$$
we have 
$$
\mathbb{E}_y \gamma^{-\tau}\le C(V(y)+1)\quad\text{for}\quad y\in Y.
$$
\end{lemma}
\noindent
{\it Proof of Lemma \ref{ol}.}\newline
Define
$$ 
\kappa=\epsilon+\rho \circ T_{\epsilon},
$$
where $T_n((y_k)_{k\in \mathbb{N}_0})=(y_{k+n})_{k\in\mathbb{N}_0}$. Fix $y\in Y$, $\alpha\in (0,1)$ and $r>1$ such that 
$(\lambda \alpha)^{-\frac{1}{r-1}}<\beta^{-1}$. The strong Markov property and the H\" older inequality for every $y\in Y$ give 
$$
\begin{aligned}
\mathbb{E}_y {\bf 1}_{\{\epsilon <\infty\}}\lambda^{-\frac{\kappa}{r}} &\le
[\mathbb{E}_{y} ({\bf 1}_{\{\epsilon <\infty\}}(\lambda\alpha)^{-\frac{\epsilon}{r}})^{\frac{r}{r-1}}]^{\frac{r-1}{r}}
[\mathbb{E}_{y} ({\bf 1}_{\{\epsilon <\infty\}}\alpha^{\epsilon}\lambda^{-\rho\circ T_{\epsilon}})]^{\frac{1}{r}}\\
&\le (\mathbb{E}_y {\bf 1}_{\{\epsilon <\infty\}}\beta^{-\epsilon})^{\frac{r-1}{r}}
[\mathbb{E}_{y} ({\bf 1}_{\{\epsilon <\infty\}}
\alpha^{\epsilon}\mathbb{E}_y(\lambda^{-\rho\circ T_{\epsilon}}|\mathcal{F}_{\epsilon}))]^{\frac{1}{r}}\\
&= (\mathbb{E}_y {\bf 1}_{\{\epsilon <\infty\}}\beta^{-\epsilon})^{\frac{r-1}{r}}
[\mathbb{E}_{y} ({\bf 1}_{\{\epsilon <\infty\}}\alpha^{\epsilon}\mathbb{E}_{Y_{\epsilon}^y}(\lambda^{-\rho}))]^{\frac{1}{r}}\\
&\le (\mathbb{E}_y {\bf 1}_{\{\epsilon <\infty\}}\beta^{-\epsilon})^{\frac{r-1}{r}}
[\mathbb{E}_{y} ({\bf 1}_{\{\epsilon <\infty\}}\alpha^{\epsilon}C_{\rho}(V(Y_{\epsilon}^y)+1))]^{\frac{1}{r}},
\end{aligned}
$$
where $\mathcal{F}_{\epsilon}$ is the $\sigma$-algebra generated by $\epsilon$. Since $\sup_{y\in A} V(y)<\infty$ and $V$ satisfies
$$
\mathbb{E}_y ({\bf 1}_{\{\epsilon <\infty\}} \alpha^{\epsilon} V(Y_{\epsilon}^y)) \le C_1 (V(y)+1)\quad\text{for}\quad y\in Y,
$$
for some $C_1>0$, taking $c=\lambda^{\frac{1}{r}}$ we obtain 
$$
\mathbb{E}_y {\bf 1}_{\{\epsilon <\infty\}} c^{-\kappa} \le C_2\quad\text{whenever}\quad y\in A,
$$
for some constant $C_2>0$. Define $\epsilon_0=0$, $\kappa_0=\rho$ and
$$
\begin{aligned}
\epsilon_n &=\kappa_{n-1} + \epsilon \circ T_{\kappa_{n-1}},\\
\kappa_n &=\kappa _{n-1}+\kappa\circ T_{\kappa_{n-1}}\quad\text{for}\quad n\ge 1.
\end{aligned}
$$
Observe that $V(Y_{\kappa_n}^y)\in A$, $Y_{\epsilon_n}^y\notin B$, $\epsilon_n\le \kappa_n\le \epsilon_{n+1}$ 
and $\kappa_n\nearrow\infty$. 
We have 
$$
\begin{aligned}
\mathbb{E}_y{\bf 1}_{\{\kappa_{n+1}<\infty\}} c^{-\kappa_{n+1}} &= 
\mathbb{E}_y[{\bf 1}_{\{\kappa_{n}<\infty\}} c^{-\kappa_n} 
\mathbb{E}_y(({\bf 1}_{\{\kappa<\infty\}}c^{-\kappa})\circ T_{\kappa_n} |\mathcal{F}_{\kappa_n})]\\
&=\mathbb{E}_y[{\bf 1}_{\{\kappa_{n}<\infty\}} c^{-\kappa_n}
\mathbb{E}_{Y_{\kappa_n}^y}({\bf 1}_{\{\kappa<\infty\}}c^{-\kappa})] \\
&\le C_2\mathbb{E}_y {\bf 1}_{\{\kappa_n<\infty\}}c^{-\kappa_n}
\end{aligned}
$$
and thus
$$
\mathbb{E}_y{\bf 1}_{\{\kappa_{n}<\infty\}} c^{-\kappa_{n}}\le C_2^{n} C_{\rho}(V(y)+1)\quad\text{for}\quad y\in Y.
$$
Define $E=Y^{\infty}\setminus (Y\times B^{\infty})$ and $B_n=\{\epsilon_n<\infty\}$. Observe that $B_{n+1}=T_{\kappa_n}^{-1} E$ 
and $B_n\in\mathcal{F}_{\epsilon_n}\subset\mathcal{F}_{\kappa_n}$. For $y\in Y$ we have
$$
\begin{aligned}
\mathbb{P}_y(B_{n+1}) &=\mathbb{P}_y(B_n\cap B_{n+1})=
\mathbb{E}_y({\bf 1}_{B_n}\mathbb{E}_{y}({\bf 1}_E\circ T_{\kappa_n}|\mathcal{F}_{\kappa_n}))\\
&=\int_{B_n}\mathbb{P}_{x_{\kappa_n}}(E)\,\mathbb{P}_y(dx)\le (1-p)\mathbb {P}_y(B_n),
\end{aligned}
$$
where $x=(x_0,x_1,...)\in Y^{\infty}$. It follows that
$$
\mathbb{P}_y(B_n)\le (1-p)^n\quad\text{for}\quad y\in Y, n\ge 1.
$$
Define 
$$
{\widehat \tau}((y_n)_{n\in\mathbb{N}_0})=\inf\{n\ge 1:\, y_n\in A,\, \forall_{k> n}\, y_k\in B\}
$$
and $D_0=\{{\widehat \tau}=\kappa_0\}$, $D_n=\{\kappa_{n-1}<{\widehat \tau}\le \kappa_n<\infty\}$, for $n\ge 1$. 
Since $B_n=\{{\widehat\tau}>\kappa_{n-1}\}$, we have $D_n\subset B_n$ for $n\ge 0$ 
and $\mathbb{P}_y({\widehat \tau}=\infty)=0$ for $y\in Y$. 
Finally, by the H\"older inequality, for $s>1$ and $y\in Y$ we obtain
$$
\begin{aligned}
\mathbb{E}_y\lambda^{-\frac{{\widehat \tau}}{s}}&\le
\sum_{n=0}^{\infty}\mathbb{E}_y({\bf 1}_{\{\kappa_n <\infty\}}\lambda^{-\frac{\kappa_n}{s}}{\bf 1}_{D_n})\\
&\le \sum_{n=0}^{\infty}[\mathbb{E}_y {\bf 1}_{\{\kappa_n <\infty\}}\lambda^{-\kappa_n}]^{\frac{1}{s}}\mathbb{P}_y(D_n)^{1-\frac{1}{s}}\\
&\le \sum_{n=0}^{\infty} [C_2^n C_{\rho}(V(y)+1)]^{\frac{1}{s}} (1-p)^{n(1-\frac{1}{s})}\\
&\le C_{\rho}^{\frac{1}{s}} (1+V(y))\sum_{n=0}^{\infty} [(\frac{C_2}{1-p})^{\frac{1}{s}}(1-p)]^n.
\end{aligned}
$$
Choosing sufficiently large s  and setting $\gamma=\lambda^{\frac{1}{s}}$ we have
$$
\mathbb{E}_y\gamma^{-{\widehat \tau}}\le C(1+V(y))\quad\text{for}\quad y\in Y.
$$
Since $\tau<{\widehat \tau}$, the proof is complete.
\null\hfill $\square$ \newline

\begin{lemma}\label{lapunov}
Let $(Y_n^y)_{n\in\mathbb{N}_0}$ with $y\in Y$ be a family of Markov chains on a metric space Y. Suppose that $V:Y\to [0,\infty)$ is a Lapunov function for their transition function $\{\pi_y:\, y\in Y\}$, i.e. there exist $a\in (0,1)$ and $b>0$ such that
$$
\int_Y V(x)\pi_y(dx)\le aV(y)+b\quad\text{for}\quad y\in Y.
$$
Then there exist $\lambda\in (0,1)$ and ${\tilde C}>0$ such that for
$$
\rho ((y_k)_{k\in{\mathbb N}_0})=\inf\{k\ge 1:\,V(y_k)<\frac{2b}{1-a}\}
$$
we have
$$
\mathbb{E}_{y} \lambda^{-\rho}\le {\tilde C}(V(y_0)+1)\quad\text{for}\quad y\in Y.
$$
\end{lemma}

\noindent
{\it Proof of Lemma \ref{lapunov}.}\newline
Chains $(Y^{y}_n)_{n\in\mathbb{N}_0}$, $y\in Y$ are defined on a common probability space $(\Omega, \mathcal{F},\mathbb{P})$. 
Fix $\alpha\in (\frac{1+a}{2},1)$ and set $V_0=\frac{b}{\alpha -a}$. Define
$$
{\tilde\rho }((y_k)_{k\in{\mathbb N}_0})=\inf\{k\ge 1 : \, V(y_k)\le V_0\}
$$
Fix $y\in Y$. Let $\mathcal{F}_n\subset\mathcal{F}$, $n\in\mathbb{N}_0$ be the filtration induced by $(Y^{y}_n)_{n\in\mathbb{N}_0}$. Define
$$
A_n=\{\omega\in\Omega:\, V(Y_i^{y}(\omega))> V_0\quad\text{for}\quad i=0,1,...,n\},\quad n\in\mathbb{N}_0.
$$
Observe that $A_{n+1}\subset A_n$ and $A_n\in\mathcal{F}_n$. By the definition of $V_0$ we have ${\bf 1}_{A_n}\mathbb{E}(V(Y_{n+1}^{y})|\mathcal{F}_n) \le {\bf 1}_{A_n}(aV(Y_n^{y})+b)<\alpha {\bf 1}_{A_n} V(Y_n^{y})$ $\mathbb{P}$-a.e. in $\Omega$. This gives
$$
\begin{aligned}
\int _{A_n} V (Y_n^{y} ) d \mathbb{P}  & \le \int _{A_{n-1}} V (Y_n^{y} )d\mathbb{P} = \int_{A_{n-1}} \mathbb{E} ( V (Y_n^{y}) | \mathcal{F}_{n-1} )d\mathbb{P} \\
& \le \int_{A_{n-1}} ( a V(Y_{n-1}^{y})  + b )d\mathbb{P} 
 \le \alpha\int_{A_{n-1}}  V(Y_{n-1}^{y}) d\mathbb{P}.
\end{aligned}
$$
By the Chebyshev inequality
$$
\begin{aligned}
 \mathbb{P}(V(Y_0^{y})> V_0,..., V(Y_n^{y})> V_0)=\int_{A_{n-1}}\mathbb{P}(V(Y_n^{y})>V_0|\mathcal{F}_{n-1})d\mathbb{P}\\
\le V_0^{-1}\int_{A_{n-1}} \mathbb{E}(V(Y_n^{y})|\mathcal{F}_{n-1})d\mathbb{P}\le \alpha^nV_0^{-1}(aV(y)+b),
\end{aligned}
$$
thus for some $C>0$ we have
$$
\mathbb{P}_{y} ({\tilde \rho} >n) \le \alpha^n C (V(y)+1),\qquad n\in\mathbb{N}_0.
$$
Fix $\gamma\in (0,1)$ and observe that for $\lambda=\alpha^{\gamma}$ we have
$$
\mathbb{E}_{y} \lambda^{-{\tilde\rho}}  \le 2+\sum\limits_{n=1}^{\infty} \mathbb{P}_{y} (\lambda^{-{\tilde\rho}}>n)
\le 2+ \frac{C(V(y)+1)}{\alpha}\sum\limits_{n=1}^{\infty} n^{-\frac{1}{\gamma}}={\tilde C}(V(y)+1)
$$
for properly chosen ${\tilde C}$. Since $\rho\le {\tilde\rho}$, the proof is finished.
\null\hfill $\square$ \newline

\noindent
{\bf Proof of Theorem \ref{t1}.}\newline
{\bf Step I:} Define new metric ${\bar d}(x,y)=d(x,y)^{\nu}$ and observe that for ${\bar D}_r=\{(x,y)\in X^2:\, {\bar d}(x,y)<r\}$ we have $D_R={\bar D}_{{\bar R}}$ with ${\bar R}=R^{\nu}$. By the Jensen inequality (\ref{A2}) takes form
\begin{equation}\label{A2prime}
\int_{ X^2} {\bar d} (u,v){\mathbf Q}_{x,y}(d u,d v) \le {\bar \alpha} {\bar d}(x,y)\qquad for \qquad (x,y)\in F,
\end{equation}
with ${\bar \alpha}=\alpha^{\nu}$. Assumption {\bf A3} implies that
\begin{equation}\label{A3prime}
1-\| {\mathbf Q}_{x,y} \| \le l {\bar d}(x,y)\ \qquad\text{and}\qquad {\mathbf Q}_{x,y} (D_{{\bar \alpha}{\bar  d}(x,y)} )\ge\delta
\end{equation}
for $(x,y)\in  F$.\newline
{\bf Step II:} Observe, that if $b\in\mathcal{M}_{fin}(X^2)$ satisfies $supp\, b\subset F$ then (\ref{A3prime}) implies
$$
\|\mathbf{Q}_b\|\ge \|b\|- l \int_{X^2} {\bar d}(u,v) b(du, dv).
$$
Iterating the above inequality we obtain
$$
\|\mathbf{Q}_b^{0,...,n}\|\ge\|b\|-\frac{l}{1-{\bar \alpha}}\int_{X^2} {\bar d}(u,v) b(du,dv).
$$
If $supp \, b\subset\{(u,v)\in X^2:\, {\bar d}(u,v)<\frac{1-{\bar \alpha}}{2l}\}\cap F$ then from (\ref{Qinfty}) it follows that 
\begin{equation}\label{1pt1}
\|\mathbf{Q}_b^{\infty}\|\ge \frac{1}{2}\|b\|.
\end{equation}
Set $R_0=\sup\{{\bar d}(x,y):\, L(x)+L(y)<R\}<\infty$ and 
$n_0=\min\{n\in\mathbb{N}_0:\,{\bar \alpha^n} R_0<\frac{1-{\bar \alpha}}{2l}\}$. Now (\ref{A3prime})
implies that for $(x,y)\in F$ such that $L(x)+L(y)<R$ we have
$$
\mathbf{Q}_{x,y}^{n_0}(\{(u,v)\in X^2:\, {\bar d}(u,v)<\frac{1-{\bar \alpha}}{2l}\}\cap F)>\delta^{n_0}
$$
and finally (\ref{1pt1}) gives
\begin{equation}\label{2pt1}
\|\mathbf{Q}_{x,y}^{\infty} \|\ge\frac{1}{2}\delta ^{n_0}.
\end{equation}
{\bf Step III:} Define $\tilde\rho ( (x_n,y_n)_{n\in\mathbb{N}_0})=\inf\{n\ge 1 :\, L(x_n)+L(y_n)<\frac{4c}{1-\lambda}\}$.
Since $L(x)+L(y)$ is a Lapunov function for a Markov chain in $X^2$ with transition probabilities 
$\{\mathbf{B}_{x,y}:\, x,y\in X\}$, Lemma \ref{lapunov} shows that there exist constants $\lambda_0\in(0,1)$ and $C_0$ such that
\begin{equation}\label{3pt1}
\mathbb{E}_{x,y}\, \lambda_0 ^{-\tilde\rho}\le C_0 (L(x)+L(y)+1)\qquad\text{for}\qquad (x,y)\in X^2.
\end{equation}
Define $A=\{(x,y,\theta)\in {\widehat X}:\, (x,y)\in F\quad\text{and}\quad L(x)+L(y)<R\}$ and
$$
\rho ((x_n,y_n,\theta_n)_{n\in\mathbb{N}_0})=\inf\{n\in\mathbb{N}_0:\, (x_n,y_n,\theta_n)\in A\}.
$$
Since $\rho\le\tilde\rho +\kappa \circ T_{\tilde\rho}$, where
$T_{\tilde\rho} ((x_n,y_n,\theta_n)_{n\in\mathbb{N}_0})=
(x_{n+\tilde\rho},y_{n+\tilde\rho},\theta_{n+\tilde\rho})_{n\in\mathbb{N}_0},
$
 an argument similar to that in the proof of Lemma \ref{ol} shows that there exist $\lambda\in (0,1)$ such that 
$$
\mathbb{E}_{x,y,\theta}\, \lambda^{-\rho}\le \tilde C C_0 (L(x)+L(y)+1)\qquad\text{for}\qquad x,y\in X,\theta\in\{0,1\}.
$$
Define $B=\{(x,y,\theta)\in{\widehat X}:\, \theta_n=1\}$ and 
$$
\epsilon((x_n,y_n,\theta_n)_{n\in\mathbb{N}_0}=\inf\{n\ge1:\, (x_n,y_n,\theta_n)\not\in B\}.
$$
From Step II we obtain $\mathbb{P}_{x,y,\theta} (B)\ge\frac{1}{2}\delta^{n_0}$ for 
$(x,y,\theta)\in A$.
From (\ref{A2prime}) and (\ref{A3prime}) it follows that
$$
\begin{aligned}
{\widehat {\mathbf B}}_{x,y,\theta}(\epsilon=n)=& \int_{{\widehat X}^n} {\widehat {\mathbf R}}_{z_{n-1}}({\widehat X})
\, {\widehat {\mathbf Q}}_{x,y,\theta}^{0,...,n-1}(dz_0,...,dz_{n-1})\\
=& ||{\mathbf Q}^{n-1}_{\delta_{(x,y)}}||-||{\widehat Q}_{{\mathbf Q}^{n-1}_{\delta_{(x,y)}}}||
\le l\int_{X^2} {\bar d}(u,v)\, {\mathbf Q}^{n-1}_{\delta_{(x,y)}}(du,dv)\\
\le & l{\bar \alpha}^{n-1}{\bar d}(x,y)<{\bar\alpha}^{n-1}lR_0,
\end{aligned}
$$
whenever $(x,y,\theta)\in A$.
Finally Lemma \ref{ol} guarantees the existence of constants $\gamma\in (0,1),\, C_1>0$ such that for 
$$
\tau ((x_n,y_n,\theta_n)_{n\in\mathbb{N}_0})=\inf\{n\ge 1:\, \forall_{k\ge n}(x_k,y_k,\theta_k)\in B\}
$$
we have
$$
\mathbb{E}_{x,y,\theta} \, \gamma^{-\tau}\le C_1(L(x)+L(y)+1) \qquad\text{for}\qquad x,y\in X,\theta\in\{0,1\}.
$$
{\bf STEP IV:} Define sets
$$G_{\frac{n}{2}}=\{t\in(X^2\times\{0,1\})^{\infty}:\, \tau (t)\le\frac{n}{2}\}$$
and
$$H_{\frac{n}{2}}=\{t\in(X^2\times\{0,1\})^{\infty}:\, \tau (t)>\frac{n}{2}\}.$$
For every $n\in\mathbb{N}$ we have
$$
\widehat{\mathbf{B}}_{x,y,\theta}^{\infty}=\widehat{\mathbf{B}}_{x,y,\theta}^{\infty}\mid_{G_{\frac{n}{2}}} +
\widehat{\mathbf{B}}_{x,y,\theta}^{\infty}\mid_{H_{\frac{n}{2}}}\qquad\text{for}\qquad x,y\in X,\theta\in\{0,1\}.
$$
Fix $\theta=1$ and $(x,y)\in X^2$. From the fact that $\Vert \cdot\Vert_{FM} \leq \Vert \cdot \Vert_{W}$ it follows that
\begin{align*} 
& \Vert P^{*n}\delta_x-P^{*n}\delta_y\Vert_{FM}=\Vert \mathbf{P}_x^n -\mathbf{P}_y^n \Vert _{FM} \\
&=\sup\limits_{f\in\mathcal{F}} |\int_{X^2} ( f(z_1)-f(z_2)) ( pr^{\#} _n \mathbf{B}_{x,y}^{\infty} ) (dz_1,dz_2)| \\
&=\sup\limits_{f\in\mathcal{F}} |\int_{X^2} ( f(z_1)-f(z_2)) ( pr_{X^2}^{\#} pr^{\#} _n \widehat{\mathbf{B}}_{x,y,\theta}^{\infty} ) (dz_1,dz_2)|\\
&\leq \sup\limits_{f\in\mathcal{W}} |\int_{X^2} ( f(z_1)-f(z_2)) (pr ^{\#}_{X^2} pr^{\#} _n ( \widehat{\mathbf{B}}_{x,y,\theta}^{\infty}\mid _{G_{\frac{n}{2}}}) ) (dz_1,dz_2)| + 
2\widehat{\mathbf{B}}_{x,y,\theta}^{\infty}(H_{\frac{n}{2}} ).
\end{align*} 
From {\bf A2} we obtain
\begin{align*}
&\sup\limits_{\mathcal{W}}|\int_{X^2}(f(z_1)-f(z_2))(pr_{X^2}^{\#}pr_n^{\#}(\widehat{\mathbf{B}}_{x,y,\theta}^{\infty}\mid _{G_{\frac{n}{2}}}))(dz_1,dz_2)| \\
&\le \int_{X^2} d(z_1,z_2) (pr_{X^2}^{\#}pr_n^{\#}(\widehat{\mathbf{B}}_{x,y,\theta}^{\infty}\mid _{G_{\frac{n}{2}}}))(dz_1,dz_2)\\
&\le\alpha^{\frac{n}{2}}\int_{X^2} d(z_1,z_2) (pr_{X^2}^{\#}pr_{\frac{n}{2}}^{\#}(\widehat{\mathbf{B}}_{x,y,\theta}^{\infty}\mid _{G_{\frac{n}{2}}}))(dz_1,dz_2)
\le \alpha^{\frac{n}{2}}R.
\end{align*}
Now Step III and the Chebyshev inequality imply that
$$
\widehat{\mathbf{B}}_{x,y,\theta}^{\infty} (H_{\frac{n}{2}})\le\gamma^{\frac{n}{2}}C_1(L(x)+L(y)+1)\qquad\
\text{for}\qquad n\in\mathbb{N}.
$$
Taking $C_2=2C_1+R$ and $q=\max\{\gamma^{\frac{n}{2}},\alpha^{\frac{n}{2}}\}$ we obtain
$$
\Vert P^{*n}\delta_x-P^{*n}\delta_y\Vert_{FM}\le \gamma^nC_1(L(x)+L(y)+1)\qquad\text{for}\qquad x,y\in X,n\in\mathbb{N},
$$
and so
\begin{equation}\label{4pt1}
\Vert P^{*n}\mu-P^{*n}\nu\Vert_{FM}\le \gamma^nC_1(\int_X L(x) \mu(dx)+\int_X L(y)\nu(dy)+1)
\end{equation}
for $\mu,\nu\in\mathcal{M}^L_1 (X)$ and $n\in\mathbb{N}$.
\newline
{\bf Step V:} Observe that Step IV and {\bf A1} give
\begin {align*}
&\Vert P^{*n}\delta_x-P^{*(n+k)}\delta_x\Vert_{FM}
\le\int_{X}\Vert P^{*n}\delta_x-P^{*n}\delta_y\Vert_{FM}P^{*k}\delta_x(dy)\\
&\le q^nC_2\int_X (L(x)+L(y))P^{*k}\delta_x(dy)
\le q^n C_3 (1+L(x)),
\end{align*}
so   $(P^{*n}\delta_x)_{n\in\mathbb{N}}$ is a Cauchy sequence for every $x\in X$. Since $\mathcal{M}_1 (X)$ equipped with the norm $\Vert\cdot\Vert_{FM}$ is complete (see \cite{EK}), assumption {\bf A0} implies the existence of an invariant measure $\mu_{*}$. 
Assumption {\bf A1} gives $\mu_{*}\in\mathcal{M}_1^L(X)$. Applying inequality {\rm(\ref{4pt1})} we obtain {\rm (\ref{mix})}. Observation that $\mathcal{M}_1^L(X)$ is dense in $\mathcal{M}_1(X)$ in the total variation norm finishes the proof.
 
\null\hfill $\square$ \newline
{\it Remark.} In steps IV and V of the above proof we follow M. Hairer (see \cite{H}).

\bigskip
\section{Random iteration of functions}

Let $(X,d)$ be a Polish space and $(\Theta,\Xi)$ a measurable space with a family $\vartheta_x\in\mathcal{M}_1(\Theta)$ of distributions on $\Theta$ indexed by $x\in X$. Space $\Theta$ serves as a set of indices for a family $\{S_\theta: \theta\in\Theta\}$ of continuous functions acting on $X$ into itself.
We assume that $(\theta,x)\mapsto S_\theta(x)$ is product measurable.
In this section we study some stochastically perturbed dynamical system $(X_n)_{n\in\mathbb N_0}$. 
Its intuitive description
is the following: if $X_0$ starts at $x_0$, then by choosing $\theta_0$ at random from $\vartheta_{x_0}$
we define $X_1=S_{\theta_0}(x_0)$. Having $X_1$ we select $\theta_1$ according to the distribution
$\vartheta_{X_1}$ and we put $X_2=S_{\theta_1}(X_1)$ and so on.
More precisely, the process $(X_n)_{n\in\mathbb N_0}$ can be written as 
$$
X_{n+1}=S_{Y_n}(X_n), 
\qquad
n=0,1,\dots,
$$
where $(Y_n)_{n\in\mathbb N_0}$ is a sequence of random elements defined on a probability
space $(\Omega,\Sigma, prob)$ with values in $\Theta$ such that
\begin{equation}\label{system}
prob\,(Y_n\in B | X_n=x)=\vartheta_x(B) 
\qquad\text{for}\qquad x\in X, B\in\Xi,
n=0,1,\dots,
\end{equation}
and $X_0:\Omega\to X$ is a given random variable.
Denoting by $\mu_n$ the probability law of $X_n$, we will give 
a recurrence relation between $\mu_{n+1}$ and $\mu_n$.
To this end fix $f\in B_b(X)$ and note that 
$$
\mathbb{E} f(X_{n+1})=\int_X f d\mu_{n+1}.
$$
By (\ref{system}) we have
$$
\int_A\vartheta_x(B)\mu_n(dx)=prob (\{Y_n\in B\}\cap \{X_n\in A\})
\qquad\text{for}\quad B\in\Xi, A\in\mathcal{B}_X,
$$
hence
$$
\mathbb{E} f(X_{n+1})=\int_\Omega f(S_{Y_n(\omega)}(X_n(\omega)) prob(d\omega)=
\int_X\int_\Theta f(S_\theta(x))\vartheta_x(d\theta)\mu_n(dx).
$$
Putting $f={\bf 1}_A$, $A\in\mathcal{B}_X$, we obtain
$
\mu_{n+1}(A)=P^*\mu_n(A),
$
where
$$
P^*\mu(A)=\int_X\int_\Theta {\bf 1}_A(S_\theta(x))\vartheta_x(d\theta)\mu(dx)
\qquad\text{for}\quad \mu\in\mathcal{M}_{fin}(X), A\in\mathcal{B}_X.
$$ 
In other words this formula defines the transition operator
for $\mu_n$.
Operator $P^*$ is adjoint of the Markov operator $P:B_b(X)\to B_b(X)$ of the form
\begin{equation}\label{operator}
Pf(x)=\int_\Theta f(S_\theta(x))\vartheta_x(d\theta).
\end{equation}
We take this formula 
as the precise formal definition of considered process.
We will show that operator (\ref{operator}) 
has a unique invariant measure, provided the
following conditions hold:
\newline
{\bf B1} {\it There exists} $\alpha \in  (0,1)$ {\it such that} 
$$
\int_{\Theta} d (S_\theta(x),S_\theta(y)){\vartheta}_{x}(d \theta) \le \alpha d(x,y)\qquad\text{for}\quad x,y\in X.
$$
{\bf B2} {\it There exists}  $\bar x\in X$ {\it such that} 
$$
c:=\sup_{x\in X}\int_{\Theta} d (S_\theta(\bar x),\bar x){\vartheta}_{x}(d \theta) <\infty.
$$
{\bf B3} 
{\it The map} $x\mapsto \vartheta_x$, $x\in X$, {\it is H\"older continuous in the total variation norm, i.e.} 
{\it there exists} $l>0$ and $\nu\in (0,1]$ {\it such that} 
$$
\Vert\vartheta_x-\vartheta_y\Vert\le l\,d(x,y)^\nu 
\qquad\text{for}\quad 
x,y\in X.
$$
{\bf B4} {\it There exists} $\delta>0$ {\it such that} 
$$
\vartheta_x\wedge\vartheta_y (\{\theta\in\Theta:d(S_\theta(x),S_\theta(y))\le\alpha d(x,y)\})>\delta\qquad\text{for}\quad x,y\in X,
$$
{\it where} $\wedge $ {\it denotes the greatest lower bound in the lattice of finite measures}.

\medskip
{\it Remark.} It is well known (see \cite{Stenflo}) that 
replacing the H\"older continuity in {\bf B3} by slightly weaker condition of the Dini continuity can lead to the lack of exponential convergence.

\begin{proposition} \label{pl}
Assume {\bf B1} -- {\bf B4}. Then operator {\rm (\ref{operator})} possesses a unique invariant measure 
$\mu_{*}\in\mathcal{M}_1^1 (X)$, which is attractive in $\mathcal{M}_1(X)$. Moreover there exist $q\in (0,1)$ and $C>0$ such that
$$
\| P^{* n}\mu -\mu_{*}\|_{FM}\le q^n C(1+\int_X d(\bar x,x)\mu (dx) )
$$
for $\mu\in\mathcal{M}_1^1 (X)$ and $n\in\mathbb{N}$.
\end{proposition} 

{\it Proof.}  Define the operator $Q$ on $B_b(X^2)$ by
$$
Q(f)(x,y)=\int_\Theta f(S_\theta(x),S_\theta(y))\vartheta_x\wedge \vartheta_y(d\theta).
$$ 
Since 
$$
||\vartheta_{x^\prime}\wedge \vartheta_{y^\prime}-
\vartheta_{x}\wedge \vartheta_{y}|| \le
2 (||\vartheta_{x^\prime}-\vartheta_{x}||
+||\vartheta_{y^\prime}-\vartheta_{y}||)
$$
it follows that
$$
\begin{aligned}
|Q(f)(x^\prime,y^\prime)&-Q(f)(x,y)| \le
\int_\Theta |f(S_\theta(x^\prime),S_\theta(y^\prime))|\,||\vartheta_{x^\prime}\wedge \vartheta_{y^\prime}-
\vartheta_{x}\wedge \vartheta_{y}||(d\theta)\\
& 
+\int_\Theta |f(S_\theta(x^\prime),S_\theta(y^\prime))-f(S_\theta(x),S_\theta(y))|
\vartheta_{x}\wedge \vartheta_{y}(d\theta)\\
&  
\le 2l\sup_{z\in X^2}|f(z)| \,(d(x,x^\prime)^\nu+d(y,y^\prime)^\nu)
\\
&
+\int_\Theta |f(S_\theta(x^\prime),S_\theta(y^\prime))-f(S_\theta(x),S_\theta(y))|
\vartheta_{x}\wedge \vartheta_{y}(d\theta),
\end{aligned}
$$
for $f\in B_b(X^2)$, $x, y\in X$. 
Consequently, we see that $Q(C_b(X^2))\subset C_b(X^2)$,
by Lebesgue's dominated convergence theorem.
Put 
$$\mathcal F=\{f\in B_b(X^2): \sup_{z\in X^2}|f(z)|\le M, Q(f)\in B_b(X^2)\},$$
where $M>0$ is fixed, and observe that the family $\mathcal F$ is closed in pointwise convergence.
Therefore $\mathcal F$ consists of all Baire functions bounded by $M$.
By virtue of \cite[Theorem 4.5.2]{L} we obtain $Q(B_b(X^2))\subset B_b(X^2)$.
In particular, 
for the family $\{Q_{x,y}: {x,y\in X}\}$ of (subprobability) measures  given by
$$
Q_{x,y}(C)=\int_\Theta {\bf 1}_C(S_\theta(x),S_\theta(y))\vartheta_x\wedge\vartheta_y(d\theta),
$$
we have that maps $(x,y)\mapsto Q_{x,y}(C)$ are measurable for every $C\in\mathcal B_{X^2}$.

Arguing similarly as above we show that (\ref{operator})
is well defined Feller operator. It has
Lapunov function $L(x)=d(x,\bar x)$, since
$$
\int_\Theta d(S_\theta(x),\bar x)\vartheta_x(d\theta)\le \alpha d(x,\bar x)+c.
$$

Now, observe that
$$
\Vert Q_{x,y}\Vert=\vartheta_x\wedge\vartheta_y(\Theta)=1-\sup_{A\in\Xi}\{\vartheta_y(A)-\vartheta_x(A)\}\ge 1-l\, d(x,y)^\nu
$$
for $x,y\in X$.
Moreover, we have
$$
\int_{X^2}d(u,v) Q_{x,y}(du,dv)=\int_\Theta d(S_\theta(x),S_\theta(y))\vartheta_x\wedge\vartheta_y(d\theta)\le \alpha d(x,y),
$$
and
$$
Q_{x,y}(D_{\alpha d(x,y)})=
\vartheta_x\wedge\vartheta_y (\{\theta\in\Theta:d(S_\theta(x),S_\theta(y))\le\alpha d(x,y)\})>\delta
$$
for $x,y\in X$.
In consequence {\bf A0} -- {\bf A3} are fulfilled.
The use of Theorem \ref{t1} (see also Remark concerning assumption {\bf A4}) ends the proof.
\null\hfill $\square$ \bigskip

\section{Perpetuities with place dependent probabilities}

Let $X=\mathbb R^d$ and $G=\mathbb R^{d\times d}\times\mathbb R^d$, 
and consider the function $S_\theta:X\to X$ defined by
$S_\theta(x)=M(\theta)x+Q(\theta)$, where 
$(M,Q)$ is a random variable on $(\Theta,\Xi)$ with values in $G$.
Then  (\ref{operator}) may be written as
\begin{equation}\label{perp1}
Pf(x)=\int_G f(mx+q) d \vartheta_x\circ (M, Q)^{-1} (m,q)
\end{equation}
This operator is  
connected with the random difference equation of the form
\begin{equation}\label{equation}
\Phi_n=M_n \Phi_{n-1}+Q_n,\qquad n=1,2,\dots,
\end{equation}
where $(M_n,Q_n)_{n\in\mathbb N}$ is a sequence of independent random variables
distributed as 
$(M,Q)$.   
Namely, the process $(\Phi_n)_{n\in\mathbb N_0}$ is homogeneous Markov chain with the transition kernel $P$ given by
\begin{equation}\label{op3}
Pf(x)=\int_G f(mx+q)d\mu(m,q),
\end{equation}
where $\mu$ stands for the distribution of $(M,Q)$.
Equation (\ref{equation}) arises in various disciplines as economics, physics, nuclear technology, biology, sociology  (see e.g. \cite{V1979}).
It is closely related to a sequence of backward iterations $(\Psi_n)_{n\in\mathbb N}$, given by
$
\sum_{k=1}^n M_1\dots M_{k-1}Q_k, n\in\mathbb N
$
(see e.g. \cite{GM2000}).
Under conditions ensuring the almost sure convergence of the sequence
$(\Psi_n)_{n\in\mathbb N}$ the limiting random variable 
\begin{equation}\label{perpetuity}
\sum_{n=1}^\infty M_1\dots M_{n-1}Q_n
\end{equation}
is often called  {\it perpetuity}. 
It turns out that the probability law of  (\ref{perpetuity})
is a unique invariant measure for (\ref{op3}).
The name perpetuity
comes from perpetual payment streams
and recently gained
some popularity in the literature on stochastic recurrence equations
(see \cite{EKM1997}).	 
In the insurance context
a perpetuity represents the present value of a permanent commitment
to make a payment at regular intervals, say annually, into the future forever.
The $Q_n$ represent annual payments, the $M_n$ cumulative discount factors. 
Many interesting examples of perpetuities can be 
found in \cite{A2009}.  
Due to significant papers \cite{K1973}, \cite{G1974}, \cite{V1979} and \cite{GM2000}
we have complete (in the dimension one) characterization
of convergence of perpetuities.
The rate of this convergence has recently been extensively studied by many authors (see for instance \cite{BJMW}-\cite{BDG}, \cite{MM}).
The main result of this section concerns the rate of convergence
of the process $(X_n)_{n\in\mathbb N_0}$ associated with the operator
$P:B_b(\mathbb R^d)\to B_b(\mathbb R^d)$ given by
\begin{equation}\label{placedependent}
Pf(x)=\int_G f(mx+q)d\mu_x(m,q),
\end{equation}
where $\{\mu_x: x\in\mathbb R^d\}$ is a family of Borel probability measures on $G$.
In contrast to $(\Phi_n)_{n\in\mathbb N_0}$,
the process $(X_n)_{n\in\mathbb N_0}$ moves by choosing at random $\theta$ from a measure depending on $x$.
Taking into considerations the concept of perpetuities we may say that $(X_n)_{n\in\mathbb N_0}$ forms {\it a perpetuity with place dependent probabilities}.

\begin{corollary} \label{p2}
Assume that 
$\{\mu_x: x\in\mathbb R^d\}$ is a family of Borel probability measures on $G$ such that$\,^1$
\begin{equation}\label{conditions}
\alpha:=\sup_{x\in\mathbb R^d}\int_G||m||d\mu_x(m,q)<1,\qquad c:=\sup_{x\in\mathbb R^d}\int_G |q|d\mu_x(m,q)<\infty.
\end{equation}
Assume moreover that 
the map $x\mapsto \mu_x$, $x\in X$, is H\"older continuous in the total variation norm
and there exists $\delta>0$ such that
$$
\mu_x\wedge\mu_y (\{(m,q)\in G:||m||\le\alpha\})>\delta\qquad\text{for}\quad
x,y\in\mathbb R^d.
$$
Then operator {\rm(\ref{placedependent})}
possesses a unique invariant measure 
$\mu_{*}\in\mathcal{M}_1^1 (\mathbb R^d)$, which is attractive in 
$\mathcal{M}_1(\mathbb R^d)$.
Moreover there exist $q\in (0,1)$ and $C>0$ such that
$$
\| P^{* n}\mu -\mu_{*}\|_{FM}\le q^n C(1+\int_{\mathbb R^d} |x|\mu (dx) )
$$
for $\mu\in\mathcal{M}_1^1 (\mathbb R^d)$ and $n\in\mathbb{N}$.
\end{corollary}

\footnotetext[1]{$||m||=\sup\{|mx|:x\in\mathbb R^d, |x|=1\}$, and $|\cdot|$ is Euclidean norm in $\mathbb R^d$}

\noindent The proof of corollary is straightforward application of Proposition \ref{pl}.
We leave the details to the reader.
We finish the paper by giving an example to illustrate Corollary \ref{p2}. 

\noindent {\bf Example.} Let $\nu_0$, $\nu_1$ be distributions on $\mathbb R^2$.
Assume that $p, q:\mathbb R\to [0,1]$ are Lipschitz functions (with Lipschitz constant
$L$) summing up to 1, and
$p(x)=1$, for $x\le 0$, 
$p(x)=0$, for $x\ge 1$.
Define $\mu_x$ by
$$
\mu_x=
p(x)\nu_0+q(x)\nu_1,\qquad x\in\mathbb R.
$$
Then:

\begin{enumerate}
\item 
$
\Vert\mu_x-\mu_y\Vert\le 2L |x-y| \qquad\text{for}\quad x,y\in\mathbb R.
$
\item 
If
$
\int_{\mathbb R^2}|m|d\nu_i(m,q)<1
$
and 
$\int_{\mathbb R^2}|q|d\nu_i(m,q)<\infty$
for $i=0,1$, then 
(\ref{conditions}) holds.
\item 
For every 
$A\in\mathcal B_{\mathbb R^2},$
$x,y\in\mathbb R$ we have:
$
\mu_x\wedge\mu_y (A)
\ge
\nu_0\wedge\nu_1(A)
=(\nu_0-\lambda^+)(A)=
(\nu_1-\lambda^-)(A)
\ge\max\{\nu_0(A),\nu_1(A)\}-\Vert\nu_0-\nu_1\Vert(A),
$ where $(\lambda^+,\lambda^-)$ is the Jordan decomposition of $\nu_1-\nu_0$.

\end{enumerate}

\end{document}